\documentclass{article}
\usepackage{amsbsy,amssymb,amscd,amsfonts,latexsym,amstext,delarray,
amsmath, diagrams}  \headsep=15pt
\def\hpi{\hat{\pi}}

\def\cX{{\cal X} }

\def\bX{ \bar{X}}

\def\dim{{ \mbox{dim} }}

\def\Hom{{ \mbox{Hom} }}

\def\ra{{ \rightarrow }}
\def\da{{ \downarrow }}

\def\hra{{ \hookrightarrow }}
\def\da{{ \downarrow }}

\def\bs{ \backslash}

\def\G{{ \Gamma }}
\def\Gal{{ \mbox{Gal} }}

\def\bQ{\bar{\Q}}

\def\cE{ {\cal E}}

\def\Z{{ \mathbb{Z}}}

\def\bq{\begin{quote}}
\def\eq{\end{quote}}

\newtheorem{thm}{Theorem}[section]

\newtheorem{claim}[thm]{Claim}

\def\Q{\mathbb{Q}}

\def\invlim{\varprojlim}

\def\be{\begin{equation}}
\def\ee{\end{equation}}

\def\cL{{\cal L}}
\def\bcL{\bar{\cal L}}
\def\bpi{\bar{\pi}}

\def\cL{ {\cal L}}

\def\s{ \sigma}

\def\bG{\bar{G}}

\def\L{\Lambda}
\def\bL{\bar{\Lambda}}
\def\bchi{\bar{\chi}}
\def\bpi{\bar{\pi}}

\title{
$p$-adic $L$-functions and Selmer varieties
associated to elliptic curves with
complex multiplication}
\author{Minhyong Kim}
\begin{document}
\maketitle
The study of non-abelian fundamental groups renders it plausible that the
principle of Birch and Swinnerton-Dyer, whereby 
non-vanishing of $L$-values, in some appropriate sense, accounts for the finiteness of integral points, can
eventually be extended to hyperbolic curves.
Here we will discuss the very simple case of a genus 1 hyperbolic curve $X/\Q$
obtained by removing the origin from an elliptic curve $E$
defined over $\Q$ with complex multiplication by an imaginary quadratic field
$K$.
Denote by $\cE$ a Weierstrass minimal model of $E$ and by
$\cX$ the integral model of $X$
obtained as the complement in $\cE$ of the closure of the origin.
Let $S$ be a set of primes including the infinite place and those of
bad-reduction for $\cE$.
We wish to examine the theorem of Siegel, asserting the finiteness of
$\cX(\Z_S)$, the $S$-integral points of $\cX$,  from the point of view of
fundamental groups and Selmer varieties. In particular,
we show how the finiteness of points can be proved using  `the method of
Coates and Wiles' which, in essence, makes use of the non-vanishing of $p$-adic
$L$-functions arising from the situation.

That is to say, in studying the set $\cE(\Z_S) (=\cE(\Z)=E(\Q))$,
Coates and Wiles showed the special case of the conjecture of Birch and
Swinnerton-Dyer by deriving the finiteness of $\cE(\Z_S)$ from the
non-vanishing of $L(E/\Q,s)$ at $s=1$. Of course the $L$-function can vanish
at 1 in general,
in which case $\cE(\Z_S)$ is supposed to be infinite. But we know that
$\cX(\Z_S)$ is always finite. From the perspective of this paper, this is
a consequence of the fact that appropriate $p$-adic $L$-functions
have only finitely many zeros.
More precisely, if we choose a prime $p$ that splits as $p=\pi \bpi$ in
$K$ and let $\cL $ and $\bcL $ denote the
$p$-adic $L$-functions associated to the $\pi$ and $\bpi$-power torsion points of $E/K$ \cite{rubin}, we know that they have only finitely many zeros.
And then, the non-vanishing of $L$-functions forces the vanishing of infinitely many  $\Q_p-$Selmer groups
 for a family of Galois representations naturally associated to   $X$.
The motivic tool used to put this information together
 in the present approach is a natural
quotient $$W$$ of the $\Q_p$-pro-unipotent fundamental
group $U$ of $X$ with base point at an $S$-integral point $b$ (which we assume to exist) and, as usual, its
further quotients $W_n$ modulo the descending central series.
Integral points of $\cX$ give rise to torsors for the $W_n$ that are classified by
a projective system of global Selmer varieties
$$H^1_f(\G,W_n)$$
where $\G=\Gal(\bQ/\Q)$ is the Galois group over $\Q$ of an algebraic closure.
That is, there is  a diagram
$$\begin{array}{ccc}
\cX(\Z_S) & \hra & \cX(\Z_p) \\
\da & & \da \\
H^1_f(\G,W_n) & \ra & H^1_f(\G_p,W_n)
\end{array}
$$
obtained from the formalism of the fundamental group
that associates to each point $x$ the $U$-torsor of paths from
$b$ to $x$ and then pushes it out to a $W_n$-torsor. Here,
$\G_p=\Gal(\bQ_p/\Q_p)$ embedded into $ \G$ as a decomposition group at $p$.

Let $G=\Gal(K(E[\pi^{\infty}])/K)$ and $\bG=\Gal(K(E[\bpi^{\infty}])/K)$,
and let $\L=\Z_p[[G]]$, $\bL=\Z_p[[\bG]]$ be the corresponding
Iwasawa algebras, so that $\cL \in \L$ and $\bcL\in \bL$.
Denote by $$\chi:\L \ra \Q_p,\ \ \bchi:\bL \ra \Q_p$$
the homomorphisms corresponding to the Galois actions
on $V_{\pi}(E)=(\invlim E[\pi^n])\otimes \Q)$
and $V_{\bpi}(E)=(\invlim E[\bpi^n])\otimes \Q)$.
Finally, let
$r=
\dim H^1_f(\G, V_p(E))$ and $s=|S|$.
\begin{thm} We have the inequality of dimensions
$$\dim H^1_f(\G,W_n) < \dim H^1_f(\G_p,W_n)$$
for all $n$ sufficiently large.
\end{thm}
The fact that the $p$-adic $L$-functions
can have at most finitely many zeros is exactly the
required input for this theorem. But
the location of the zeros, of course, is a highly non-trivial
issue. On the other hand,
\begin{thm} Suppose $\chi^k(\cL)\neq 0$ and $\bchi^k(\bcL)\neq 0$
for each $k<0$. Then
$$\dim H^1_f(\G,W_n) < \dim H^1_f(\G_p,W_n)$$
for all $n\geq r+s+1$.
\end{thm}
I am informed by John Coates that the non-vanishing in the hypothesis
is a conjecture appearing in folklore.

The way the implication works out is that each non-vanishing
of an $L$-value implies the vanishing of some $H^2$ in Galois
cohomology, and $W$ is constructed so as to avoid the
complications  that arise when working directly with
$U$.
The finiteness of global points then follows in a straightforward
way as in previous work (e.g \cite{kim2}) whereby the inequality implies the
existence of certain $p$-adic iterated integrals
that vanish on the global points. Since we hope the method
will eventually lead to a direct construction of
a $p$-adic analytic function
that annihilates the global points, the more refined statement
of the second theorem seems worth making explicit.
Of course  in the present work the
main emphasis is the sequence of implications
\bq
Non-vanishing of $L$--values $\Rightarrow$ control of Selmer varieties
$\Rightarrow$ finiteness of global points,
\eq entirely parallel to the case of elliptic curves,
with just the replacement of Selmer groups by Selmer
varieties.

\medskip
A word of caution regarding the notation: At the urging of Richard Hain,
the indexing of the finite-dimensional quotients
of $U$ has been shifted. So our $U_n$ is $U_{n+1}$ from the papers \cite{kim2} and \cite{KT},
for example. However, the scheme here is consistent with that of
\cite{kim1}.

\section{A quotient of the unipotent fundamental group}
The quotient in question is  constructed as follows.

Let $U=\pi_1^{un}(\bX,b)$ be the $\Q_p$-pro-unipotent completion of
$\hpi_1(\bX,b)$ (see \cite{deligne}) and let $U^n$ denote the descending central
series, normalized so that $U^1=U$.
Define $U_n=U^{n+1}\bs U$. We then have exact sequences
$$0\ra U^{n+1}\bs U^n \ra U_n \ra U_{n-1}\ra 0$$
for $n\geq 1$.
Denote by $L$ the Lie algebra of $U$ with descending
central series $L^n$. Thus, we have  natural isomorphisms
$$U^{n+1}\bs U^n \simeq L^{n+1}\bs L^n$$
compatible with the action of $\G$.
Since $\hpi_1(\bX,b)$ is  pro-finite free on two generators, $L$
is the pro-nilpotent completion of the free Lie algebra on
two generators, where the generators can be any two elements
projecting to a basis of
$L_1=H_1(\bX,\Q_p)$.
Therefore, $L$ comes with a natural grading (not compatible with
the Galois action)
$$L=\overline{\oplus_{n=1}^{\infty} L(n)}$$
where $L(n)$ is generated by the Lie monomials of
degree $n$ in the generators \cite{serre}.
On the other, in the current situation,
we have $$L_1=V_p(E)\simeq V_{\pi}(E)\oplus V_{\bpi}(E)$$
so that
$L$ even has a bi-grading
$$L=\overline{\oplus L_{i,j}}$$
That is,
taking  $e$ and $f$ to be elements in $L_1$ that map to bases of
$V_{\pi}$ and $V_{\bpi}$ respectively,
$L_{i,j}$ is spanned by Lie monomials that have $i$ number of $e$'s
and $j$ number of $f$'s, e.g.,
$$ad^{i-1}(e)(ad^{j-1}(f)([e,f]))$$
This bi-grading also induces a filtration
$$L_{\geq n,\geq m}:=\overline{\oplus_{i\geq n,j\geq m}L_{i,j}}$$
by Lie ideals.
Let $N=\Gal(\bQ/K)\subset \G$ so that
$\G=N\s$, where $\s$ is complex conjugation.
Then if $x\in N$, we have
$$xe=\chi(x)e+z$$ and $$xf=\bchi(x)f+z'$$ for $z,z'\in L^2$.
But $L^2\subset L_{\geq 1,\geq 1}$. Hence, if
$l\in L_{i,j}$, then an easy induction shows that
$$xl=\chi(x)^i\bchi(x)^jl+z$$ for
$z\in L_{\geq i+1,\geq j+1}$. In particular, each
$L_{\geq n,\geq m}$ is stabilized by $N$. Similarly,
$$\s(L_{\geq i,\geq j}) \in L_{\geq j,\geq i}$$so we see that
$$
L_{\geq n,\geq n}$$
is stabilized by $\G$ for each $n$.
 Therefore, we get a a quotient Lie algebra
$$L \ra {\cal W}:=L/L_{\geq 2,\geq 2}\ra 0$$
and a corresponding quotient group
$$U \ra W \ra 0$$
with a compatible $\G$-action.
If we choose the ordering $e<f$, then
$[e,f]$ is a Hall basis \cite{serre} for $L^3\bs L^2$, and
inside the Hall basis  for $L^{n+1}\bs L^n$, $n\geq 3$,
are the elements
$$ad^{n-2}(e)([e,f])$$ and $$ad^{n-2}(f)([e,f]).$$
All the other basis elements are clearly in
$L_{\geq 2,\geq 2}$. Thus,
$$W^{n+1}\bs W^n\simeq {\cal W}^{n+1}\bs {\cal W}^n $$
is generated by the class of
$ad^{n-2}(e)([e,f])$ and $ad^{n-2}(f)([e,f])$.
That is, as $\G$-modules, we have
$$W^{n+1}\bs W^n\simeq \Q_p(\chi^{n-2}(1))\oplus \Q_p(\bchi^{n-2}(1))$$
for $n\geq 2$,
where $\s$ acts by sending the generator
$ad^{n-2}(e)([e,f])$ to $-ad^{n-2}(f)([e,f])$.
\section{The unipotent Albanese map}
There are various ways to see that
$W$ is unramified outside $S$ and crystalline at $p$.
For example, by construction, the coordinate ring of $W$
is a sub-ring of that of $U$, and hence, the conditions
of being unramified or crystalline (\cite{kim2}, section 2) are inherited from $U$.

Now we examine the unipotent Albanese map \cite{kim2}
$$\cX(\Z_S) \ra H^1(\G,U_n),$$ obtained by
associating to $x\in \cX(\Z_S)$ the class of the $U_n$-torsor
$$P(x):=\pi_1^{un}(\bX;b,x)$$ of unipotent paths from
$b$ to $x$. We will continue this map
to
$$\cX(\Z_S)\ra H^1(\G,U_n)\ra H^1(\G,W_n)$$
obtained by composing with the quotient map.
At the level of torsors, the map
$$ H^1(\G,U_n)\ra H^1(\G,W_n)$$
is given by the pushout:
$$Z\mapsto (Z\times W_n)/U_n$$
Since the condition of being crystalline at $p$ (\cite{kim2}, section 2)
merely signifies that a torsor has a
$\G_p$-invariant $B_{cr}$ point,
this condition is clearly preserved by push-out.
Thus, we have an induced map
$$H^1_f(\G_p, U_n)\ra H^1_f(\G_p, W_n)$$
where the subscript $f$ refers exactly to the
subset
$$
H^1_f(\G_p, W_n)\subset H^1 (\G_p, W_n)
$$
of cohomology classes that are crystalline.
Similarly, the condition of beling
unramified at $v\notin T=S\cup\{p\}$
is preserved under pushout. Let $\G_T$
denote the Galois group of the maximal
extension of $\Q$ unramified outside $T$.
By \cite{kim1}, the system
$$H^1(\G_T, W_n)$$
has the structure of a pro-algebraic variety over $\Q_p$,
as does the sub-system
$$H^1_f(\G_T, W_n)\subset H^1(\G_T, W_n)$$
classifying torsors that are unramified outside $T$
and crystalline at $p$.
Thus, we also have
an induced map of global Selmer varieties
$$H^1_f(\G, U_n)\ra H^1_f(\G, W_n)$$
giving rise to a commutative diagram
$$\begin{array}{ccc}
\cX(\Z_S)&\ra &\cX(\Z_p)\\
\da & & \da \\
H^1_f(\G, W_n)&\ra &H^1_f(\G_p, W_n)
\end{array}$$
Since each map
$$\cX(\Z_p) \ra H^1_f(\G_p,U_n)$$
has Zariski dense image, so do the
maps
$$\cX(\Z_p) \ra H^1_f(\G_p,W_n).$$
As in \cite{kim1} and \cite{kim2}, this denseness is
an important ingredient in extracting Diophantine
finiteness out of the theorems.
\section{Proof of the theorems}
The proof is now straightforward. Recall that we have
a  sequence
$$0\ra H^1(\G_T, W^{n+1}\bs W^n) \ra H^1(\G_T, W_n) \ra H^1(\G_T, W_{n-1})$$
which is exact in that the vector group
kernel acts on the middle term with quotient variety being the image of the second map.
Furthermore, if we examine
$$H^1_f(\G_p,  W^{n+1}\bs W^n)\simeq H^1_f(\G_p,  \Q_p(\chi^{n-2}(1))\oplus \Q_p(\bchi^{n-2}(1)))$$
we see that all classes are crystalline for $n\geq 3$. This is
because, for example,  the crystalline classes in
$H^1(\G_p, \Q_p(\chi^{n-2}(1)))$ are classified by
$$D^{DR}(\Q_p(\chi^{n-2}(1)))/F^0$$
where $D^{DR}(\cdot)=[(\cdot)\otimes B_{DR}]^{\G_p}$
is Fontaine's Dieudonn\'e functor \cite{BK}.
But $$D^{DR}(\Q_p(\chi^{n-2}(1)))=D^{DR}(\Q_p(\chi^{n-2}))(1)$$ and
the Hodge filtration on $D^{DR}(\Q_p(\chi^{n-2}))$ is non-positive.
Therefore, the Hodge filtration on 
$D^{DR}(\Q_p(\chi^{n-2}(1)))$ is strictly negative.
 Then a simple dimension count shows that
 $$H^1_f(\G_p, \Q_p(\chi^{n-2}(1)))=H^1(\G_p, \Q_p(\chi^{n-2}(1)))$$
 for $n\geq 3$ and the previous exact sequence can be re-written
$$0\ra H^1(\G_T, W^{n+1}\bs W^n) \ra H^1_f(\G_T, W_n) \ra H^1_f(\G_T, W_{n-1})$$
for $n\geq 3$. For $n=2$, we have
$$0\ra H^1_f(\G_T, \Q_p(1)) \ra H^1_f(\G_T, W_2) \ra H^1_f(\G_T, W_1)$$
and
$$H^1_f(\G_T, \Q_p(1))\simeq (\Z_S^*)\otimes \Q_p$$
(as in  \cite{KT}, section 2) so that
$$\dim H^1_f(\G_T, W_2)  \leq r+s-1$$
If we put this together, we get
$$\dim H^1_f(\G_T, W_n) \leq r+s-1+\Sigma_{i=3}^n \dim H^1(\G_T, W^{i+1}\bs W^i)$$
for $n\geq 3$.
As for the dimensions of the intervening $H^1$'s, we have the Euler characteristic
formula
$$\dim H^1(\G_T, W^{i+1}\bs W^i)=\dim H^2(\G_T, W^{i+1}\bs W^i)+\dim (W^{i+1}\bs W^i)^{-}$$
where the superscript refers to the subspace where $\s$ acts as $(-1)$.
But $\s$ exchanges the one-dimensional factors of
$\Q_p(\chi^{n-2}(1))\oplus \Q_p(\bchi^{n-2}(1))$ so that
$\dim (W^{i+1}\bs W^i)^{-}=1$.

Meanwhile, for the local cohomologies,
we can calculate the dimensions explicitly. Firstly, we know
that
$$\dim H^1_f(\G_p, W_2)=\dim H^1_f(\G_p, U_2)=2$$
(\cite{kim2}, section 4)
On the other hand,
$H^2(\G_p, W^{n+1}\bs W^n)=0$
for $n\geq 3$ so that the map
$$H^1(\G_p, W_{n+1}) \ra H^1(\G_p, W_{n})$$
is surjective for $n\geq 3$. As remarked above, we
also have
$H^1_f(\G_p, W^{n+1}\bs W^n)=H^1(\G_p, W^{n+1}\bs W^n)$
for $n\geq 3$. This implies that we have an exact sequence
$$0\ra H^1_f(\G_p, W^{n+1}\bs W^n) \ra H^1_f(\G_p, W_n) \ra H^1_f(\G_p, W_{n-1}) \ra 0$$
for $n\geq 3$, where each
$H^1_f(\G_p, W^{n+1}\bs W^n)$
has dimension 2. So
$$\dim H^1_f(\G_p, W_n)=2(n-2)+2=2n-2$$
for $n\geq 2$.

It remains to prove the
\begin{claim}
$$H^2(\G_T, W^{n+1}\bs W^n)=0$$
for $n$ sufficiently large.
\end{claim}
and
\begin{claim}
If $\chi^k(\L)\neq 0$ and $\bchi^k(\bL)\neq 0$
for $k<0$, then
$$H^2(\G_T, W^{n+1}\bs W^n)=0$$
for $n\geq 3$.
\end{claim}

Since 3.1 implies that
$$\dim H^1(\G_T, W^{n+1}\bs W^n)=1$$
for $n$ sufficiently large, we see that there is a constant such that
$$ H^1_f(\G_T, W_n)=C+n$$
while the local dimensions grow like $2n$. Hence we get the statement
of theorem 0.1.

On the other hand, with 3.2, we get
$$\dim H^1_f(\G_T, W_n)\leq r+s+n-2$$
so that we get the desired inequality of dimensions as soon as
$$r+s+n-2< 2n-2$$
or $n\geq r+s+1$.

We proceed to prove the claims. Clearly it suffices to
consider the Galois cohomology of
$N_T\subset G_T$, where $N_T$ is the
Galois group of the maximal extension of $K$ unramified
outside the primes dividing $T$.
For any continuous representation $M$ of
$N_T$, we define the kernel $Sha^i$ of the localization maps
on cohomology as
$$0\ra Sha^i(M) \ra H^1(N_T, M) \ra \oplus_{v|T} H^1(N_v, M)$$
where
$N_v\subset N_T$ is a decomposition group for the prime $v$.
So we have
$$0\ra Sha^2(W^{n+1}\bs W^n) \ra H^2(\G_T, W^{n+1}\bs W^n) \ra \oplus_{v| T} H^2(N_v, W^{n+1}\bs W^n)$$
By local duality,
$$H^2(N_v, W^{n+1}\bs W^n)=H^2(N_v, \Q_p(\chi^{n-2}(1))\oplus \Q_p(\bchi^{n-2}(1)))$$
$$\simeq H^0(N_v, \Q_p(\chi^{2-n})\oplus \Q_p(\bchi^{2-n}))^*=0$$
for $n\geq 3$ from which we get
$$Sha^2(W^{n+1}\bs W^n) \simeq H^2(N_T, W^{n+1}\bs W^n).$$
By Poitou-Tate duality, we have
$$Sha^2(W^{n+1}\bs W^n)\simeq Sha^1((W^{n+1}\bs W^n)^*(1))^*
\simeq Sha^1(\Q_p(\chi^{2-n}))^*\oplus Sha^1(\Q_p(\bchi^{2-n}))^*$$
But using the inflation-restriction exact sequence, we get
$$Sha^1(\Q_p(\chi^{2-n}))\simeq \Hom_{\L}(A\otimes \Q, \Q_p(\chi^{2-n}))$$
where $A$ is the Galois group of the maximal
unramified pro-$p$ extension of $K(E[\pi^{\infty}])$
split above the primes dividing $T$. Now, by \cite{rubin}, we know that the $p$-adic
$L$-function
$\cL$ annihilates $A\otimes \Q$. On the other hand,
$\L$ acts on $\Q_p(\chi^{2-n})$ through the character
$\chi^{2-n}$ and we know $\chi^{2-n}(\cL)\neq 0$ for
$n$ sufficiently large. Therefore,
$$\Hom_{\L}(A\otimes \Q, \Q_p(\chi^{2-n}))=0$$
for $n>>0$. There is a parallel argument for
 $Sha^1(\Q_p(\bchi^{2-n}))$ which then yields
 Claim 3.1. For Claim 3.2, the argument is exactly the
 same, except that the refined non-vanishing hypothesis
 implies
 $$Sha^1(\Q_p(\chi^{2-n}))\oplus Sha^1(\Q_p(\bchi^{2-n}))=0$$
 for $n\geq 3$.
 \medskip
 
 {\bf Acknowledgements:} It is a pleasure to express my gratitude
 to Kazuya Kato, Shinichi Mochizuki, Akio Tamagawa, Takeshi Tsuji
 and, especially, John Coates, conversations with whom
 were crucial to the formation of this paper.

{\footnotesize Department of Mathematics, University College London,
Gower Street, London, WC1E 6BT, United Kingdom}

\end{document}